\documentclass[12pt]{amsart}
\usepackage{amssymb}
\usepackage{amsthm}
\usepackage{amsmath}
\usepackage{amsfonts}
\newtheorem{theorem}{Theorem}[section]
\newtheorem{corollary}[theorem]{Corollary}

\newtheorem{lemma}[theorem]{Lemma}
\newtheorem{proposition}[theorem]{Proposition}

\theoremstyle{definition}
\newtheorem{example}[theorem]{Example}
\newtheorem{question}[theorem]{Question}
\newtheorem{remark}[theorem]{Remark}

\def\ppp{{\mathbb{P}}}

\def\rrr{\mathbb{R}}
\def\ccc{\mathbb{C}}
\def\zzz{\mathbb{Z}}

\def\Div{\mathrm{Div}}
\def\Pic{\mathrm{Pic}}
\def\pic{\mathrm{Pic}}
\def\Alg{\mathrm{Alg}}
\def\Num{\mathrm{Num}}
\def\ind{\mathrm{Ind}}
\def\Prin{\mathrm{Prin}}
\def\Hom{\mathrm{Hom}}
\def\deg{\mathrm{deg}}
\def\Rat{\mathrm{Rat}}

\def\Stab{\mathrm{Stab}}

\begin{document}
\title{A Question about $\pic(X)$ as a $G$-module}
\author{Daniel Goldstein}
\address{IDA/CCR, 4320 Westerra Court, San Diego, CA 92121}
\email{dgoldste@ccrwest.org}
\author{Robert M. Guralnick}
\address{Department of Mathematics, University of Southern California,
Los Angeles, CA 90089-2532}
\email{guralnic@usc.edu}
\author{David Joyner}
\address{Department of Mathematics, United States Naval Academy, 
Annapolis, MD 21402}
\email{wdj@usna.edu}
\date{\today}
\thanks{The authors are thankful Amy Ksir, David Saltman, Joe Wetherell,
and David Zelinsky for stimulating discussions.  We thank Bernard K\"och
for the reference to \cite{Lo}.}
\thanks{The second author gratefully acknowledge the support of the NSF 
grant  DMS-0140578.}

\begin{abstract}
Let $G$ be a finite group acting faithfully on
an irreducible non-singular projective curve 
defined over an algebraically
closed field $F$.
Does every $G$-invariant divisor class 
contain a $G$-invariant divisor?
The answer depends only on $G$ and not on the curve.
We answer the same question for degree $0$ divisor (classes).
We investigate the question for cycles on varieties.
\end{abstract}
\maketitle
\tableofcontents

\section{Introduction}

This paper addresses the following classical question:
Let $X$ be an irreducible non-singular projective
variety  of dimension $n$ defined over an algebraically
closed field $F$
and let $G$ be a finite subgroup of the
geometric automorphism group of $X$.
That is, $G$ is a finite group of automorphisms of 
the function field $F(X)$ that fixes $F$.
Let $D$ be an $r$-cycle and assume that its 
equivalence class $[D]$ is $G$-equivariant.
Is there always a $D'\in [D]$ which is
$G$-equivariant?
This was  addressed by Lonsted \cite{Lo}, who answered the 
question in the cyclic case for rational equivalence
but left the general case open. 
We settle the case of arbitrary finite groups 
in the case of curves and
generalize his result for varieties.
Let 
$Z^r(X)$ denote the group of $(n-r)$-cycles, let
$\Rat^r(X)$ denote the subgroup of those rationally
equivalent to $0$, let
$\Alg^r(X)$ denote the subgroup of those algebraically
equivalent to $0$, let
$\Hom^r(X)$ denote the subgroup of those homologically
equivalent to $0$, and let
$\Num^r(X)$ denote the subgroup of those numerically
equivalent to $0$. 
Let $D\in Z^r(X)$ be an $(n-r)$-cycle and assume that its 
class $[D]$ (with respect to one of these equivalences)
is $G$-equivariant.

\begin{question}Is there always a $D'\in [D]$ which is
$G$-equivariant?
\end{question}

Section 1 addresses curves. Section 2 modifies the
argument in the curve case, when possible, to varieties.
For varieties, there is an analogous question for 
algebraic equivalence classes of divisors,
and also for numerical equivalence classes of divisors.
One surprising result is that, for K3 surfaces,
if the Schur multiplier is trivial then the answer is 
yes, no matter which notion of equivalence one uses.
\section{Cohomology}
The basic idea is to use group cohomology to attack this
question. For background on cohomology,
we reference Serre \cite{S}, ch. VII. Though this is 
Lonsted's method as well, this paper was almost completely
written before this was known.
Consider the short exact sequences
\[
1\rightarrow F^\times \rightarrow F(X)^\times 
\rightarrow \Prin(X) \rightarrow 0,
\]
and
\[
0\rightarrow \Prin(X) \rightarrow \Div(X) 
\rightarrow \Pic(X) \rightarrow 0,
\]
where the additive group $\Prin(X)$ of principal divisors
may be identified with
the multiplicative group $F(X)^\times /F^\times$ via the divisor map.
Each of these is a $\zzz[G]$-module.
The covariant functor of $G$-invariants,
$M\longmapsto H^0(G,M)=M^G$ is left exact.
Therefore we have
\begin{equation}
\label{eqn:longer}
\begin{split}
1\rightarrow H^0(G,F^\times) \rightarrow H^0(G,F(X)^\times )
\rightarrow H^0(G,\Prin(X))\\ \rightarrow 
H^1(G,F^\times) \rightarrow 
H^1(G,F(X)^\times )
\rightarrow  
H^1(G,\Prin(X)) \rightarrow \\
H^2(G,F^\times)\rightarrow 
H^2(G,F(X)^\times)\rightarrow...\ ,
\end{split}
\end{equation}
and
\begin{equation}
\label{eqn:long}
\begin{split}
0\rightarrow H^0(G,\Prin(X)) \rightarrow H^0(G,\Div(X)) 
\rightarrow H^0(G,\Pic(X)) \\
\rightarrow 
H^1(G,\Prin(X)) 
\rightarrow H^1(G,\Div(X)) 
\rightarrow H^1(G,\Pic(X)) \rightarrow  ...\ .
\end{split}
\end{equation}
\section{Curves}
Let $X$ be a curve.
We claim that the answer to the above question 
is ``no'' in general. 

\begin{remark}
The answer is ``yes'' in the special case when 
$[K]$ is the canonical class.
The canonical class is $[df]$,
where $f$ is any function on the curve such that $df\ne0$. 
Clearly, there exists such an $f$ that is $G$-invariant.
\end{remark}
\begin{remark}There is a counterexample to the analogous
question for number fields. 
\end{remark}

\begin{proposition}
\label{prop:h1div=0} 
$H^1(G,\Div(X))=0$.
\end{proposition}
\begin{proof}
For each
$P\in X$, let $G_P=\{g\in G\ |\ gP=P\}$ denote its 
stabilizer (or inertia group). 
Set $L=L_P=\oplus_{g\in G/G_P} \zzz [gP].$
As an abelian group, $\Div(X)$ can be decomposed into a direct sum
of subgroups $L_P$, with one representative $P$ from each orbit
for the action of $G$ on $X$.
Note $L \cong \ind_H^G(\zzz)$, where $H\cong G_P$ is the stabilizer of one
of the basis elements. Here $G$ acts on the induced module
\[
\ind_H^G(\zzz)=\{f:G\rightarrow \zzz \ |\ f(hg)=hf(g),\ \ 
\forall h\in H,\ g\in G\},
\]
which are just the $\zzz$-valued functions on $H\backslash G$,
by right multiplication.
By Shapiro's Lemma, $H^1(G,\ind_H^G(\zzz)) \cong H^1(H,\zzz)$.
Now $H^1(H,Z)=\Hom(H,\zzz)=0$, since 
$H$ is finite.
\end{proof}

We want to know if $H^1(G,\Prin(X))$ is trivial or not, as then the
answer to the question will follow from (\ref{eqn:long}).
Observe that, due to (\ref{eqn:longer}), $H^1(G,\Prin(X)) =0$ if and only
if the map $H^2(G,F^\times)\rightarrow H^2(G,F(X)^\times)$
is injective. 
The map 
$H^2(G,F^\times)\rightarrow H^2(G,F(X)^\times)$
is defined as follows: take the cocycle defining the
extension of $G$ by $F^\times$ associated to an 
element $\alpha$ of $H^2(G,F^\times )$, use it to define an extension
of $G$ by $F(X)^\times$ in the obvious way, then
let the image of $\alpha$ be the class associated to
this extension. More precisely:
The group $H^2(G,A)$ classifies extensions $E$ of the form
\[
1\rightarrow A \rightarrow E \rightarrow G \rightarrow 1.
\]
We think of such a group $E$ as a set of pairs
$(g,a)$ with group multiplication $(g,a)(g'a')=(gg',\beta(g,g')aa')$,
for $g,g'\in G$, $a,a'\in A$ and the cocycle $\beta:G\times G\rightarrow A$
represents the associated class in $H^2(G,A)$.
So, any extension $E$ of $G$ by $F^\times$ is associated to a
cocycle $\beta:G\times G\rightarrow F^\times$.
This may be ``extended'' (apologies for the over-use
of this word) to an extension $E'$ of $G$ by $F(X)^\times$ associated
to the {\it same} cocycle. 
The map $E\longmapsto E'$ defines the map
$H^2(G,F^\times)\rightarrow H^2(G,F(X)^\times)$ in the
above long exact sequence.
In particular, the answer to our question is 
yes if and only if each non-split (central) extension of
$G$ by $F^\times$ remains non-split when 
``extended'' to $F(X)^\times$. 

Now $H^2(G,F(X)^\times))=1$ by Tsen's theorem 
(a function field 
over an algebraically closed field is a $C^1$ field).
This follows from the Corollaries on pages 96 and 
109 of Shatz \cite{Sh}. See also \S 4 and \S 7 of chapter X
in \cite{S}.
By Tsen's theorem and (\ref{eqn:longer}), 
$H^1(G,\Prin(X)) \rightarrow H^2(G,F^\times)$
is surjective. To see that the map is also
injective, note that $H^1(G,F(X)^\times )=1$,
by Hilbert's Theorem 90. So we have noted that:
\begin{lemma} 
\label{lemma:2}
$H^1(G,\Prin(X)) \cong H^2(G,F^\times)$.
\end{lemma}
When $F=\ccc$, these two lemmas imply that there is 
a $G$-equivariant representative
of every  $G$-equivariant divisor class if and only if the
Schur multiplier of $G$ is trivial.
\begin{theorem}  
\label{thrm:main}The sequence
$$\Div(X)^G \rightarrow \Pic(X)^G \rightarrow H^2(G,F^\times)\rightarrow
0$$
is exact.
In particular, there always a $D'\in [D]$ which is $G$-equivariant
if and only if $H^2(G,F^\times)=1$.
\end{theorem}

So there are many examples where the map on fixed points
fails to be surjective.

\begin{remark}
Let $p$ be the characteristic of $F.$
If $p=0$ then  $H^2(G,F^\times)$ is the Schur multiplier of $G$.
If $p>0$ then $H^2(G,F^\times)$ is the $p'$-part of the 
Schur multiplier of $G.$
\end{remark}

We give an easy corollary that has no cohomology in the statement.
\begin{corollary}
Let $G$ be a finite group. 
\begin{enumerate}
\item Let $G$ act on the curve $X$ over the complex numbers.
Suppose that $\Div(X)^G$ surjects onto $\Pic(X)^G$. 
Then $\Div(Y)^G$ surjects onto $\Pic(Y)^G$ for any curve $Y$
over any algebraically closed field.
\item Let $G$ act on the curve $X_i$ over the
algebraically closed field $K_i$ for $i=1,2$.
Suppose that $\Div(X_i)^G$ surjects onto $\Pic(X_i)^G$ for $i=1,2$,
and suppose that the characteristics of $K_1$ and $K_2$ are distinct.
Then $\Div(Y)^G$ surjects onto $\Pic(Y)^G$
surjects onto $\Pic(Y)^G$ for any curve $Y$
over any algebraically closed field.
\end{enumerate}
\end{corollary}

If the Sylow $\ell$-subgroup of $G$ is cyclic, then
the $\ell$ part of $H^2(G,F^\times)$ is trivial
(since the restriction map to the Sylow $\ell$-subgroup
 is an injection on the $\ell$ part of cohomology),
if all Sylow subgroups are cyclic, then $H^2(G,F^\times)$
is trivial.  Such groups are well known to be metacyclic.
This yields:
\begin{corollary} If 
every $\ell$-Sylow subgroup of $G$ is cyclic 
(for every prime $\ell$ dividing $|G|$),
then $\Pic(X)^G/\Div(X)^G = 0$. In particular, if $G$ is cyclic 
then for each $G$-invariant divisor class $[D]$
there is always a $D'\in [D]$ which is $G$-equivariant.
\end{corollary}
 We next want
to consider $\Pic_0(X)$ (i.e. the Jacobian) and show that the
map $\phi:\Div_0(X)^G\rightarrow \Pic_0(X)^G$ may also
fail to be surjective.
Consider the degree map $\deg$ on $\Div(X)$ and also on $\Pic(X)$.
Let $B$ be $\deg(\Div(X)^G)$.  We identify $B$ with a subgroup
of $\Pic(X)/\Pic_0(X)=\zzz$.  Clearly, $B$ contains $|G|\zzz$ and may
be bigger.  
There is no analog of Proposition \ref{prop:h1div=0}.
Instead, we have the following result.
\begin{lemma}  Let $b={\rm gcd}\, \{|G|/|I|\}$ as
$I\subseteq G$ ranges over the inertia subgroups of $G$, and set
$B=b\zzz$. Then 
\label{lemma:guralnik}
$H^1(G,\Div_0(X)) \cong \zzz/B$
\end{lemma}
\begin{proof} 
Consider the sequence
$ 0 \rightarrow \Div_0(X) \rightarrow \Div(X) \rightarrow \zzz \rightarrow 0$.
The map from $\Div(X)$ to $\zzz$ is $\deg$.  Using Proposition
\ref{prop:h1div=0} 
yields:
\[
(0\rightarrow \Div_0(X)^G\rightarrow )
\Div(X)^G \rightarrow \zzz \rightarrow H^1(G,\Div_0(X)) \rightarrow 0,
\]
as asserted.
Now we prove the claim about $b$.
Observe that a generating set of $\Div(X)^G$ are the
sum of points in a single $G$-orbit.  If $I$ is the inertia
group a point in the orbit, then the degree of this divisor
is $[G:I]$.  Thus the image of $\deg(\Div(X)^G)=B$,
where $B$ is described in the theorem.
\end{proof} 

Now consider the short exact sequence
\[
0 \rightarrow \Prin(X) \rightarrow 
\Div_0(X) \rightarrow \Pic_0(X) 
    \rightarrow 0.
\]
Taking fixed points leads to the long exact sequence for cohomology:
\[
\begin{split}
0 \rightarrow \Prin(X)^G\rightarrow 
\Div_0(X)^G \stackrel{\phi}{\rightarrow} 
\Pic_0(X)^G  \rightarrow \\
H^1(G,\Prin(X)) \rightarrow H^1(G,\Div_0(X)) \rightarrow H^1(G,\Pic_0(X)).
\end{split}
\]
This and Lemma \ref{lemma:2} proves the following result.
\begin{theorem} There is an isomorphism
\[ 
\Pic_0(X)^G/\phi(\Div_0(X)^G)  \cong
\ker \{ H^1(G,\Prin(X)) \rightarrow H^1(G,\Div_0(X))\}.
\]
In particular, if $H^2(G,F^\times)=1$ then
for each $G$-invariant divisor class $[D]$ of degree $0$
there is always a degree $0$ divisor $D'\in [D]$ which is $G$-equivariant. 
\end{theorem}
Next we want to identify the image of $H^1(G,\Prin(X))$
in $\zzz/B$, where $B$ is defined above.
Now $\deg$ induces a map from $\Pic(X)^G$ to $\zzz$.  Denote
this image by $A$ (and note that it contains $B$ the image
of $\Div(X)^G$).  Indeed, since $\Pic_0(X)^G$ is the kernel
of $\deg$, this shows that $\Pic(X)^G/(\Pic_0(X)^GI) \cong A/B$.
Thus,
$$
0 \rightarrow \Pic_0(X)^G/I_0 \rightarrow \Pic(X)^G/I \rightarrow A/B 
\rightarrow 0.
$$
This implies that the kernel of the map going from $H^1(G,\Prin(X))$ to 
$H^1(G,\Div_0(X))$ has order $|A/B|$ and since the image is cyclic,
this yields:

\begin{proposition} The image of  $H^1(G,\Prin(X))$
in $H^1(G,\Div_0(X))$ is isomorphic to $A/B$ where $A=\deg(\Pic(X)^G)$
and $B=\deg(\Div(X)^G)$.
\end{proposition}

Now $H^1(G,\Prin(X)) \cong H^2(G,F^\times)$ 
can be essentially any abelian group one wishes
(by taking direct products for example).  On the other
hand, we saw that $H^1(G,\Div_0(X))$ is a finite cyclic
group.  So there are certainly examples where this kernel
is nontrivial (and as large as one wishes). 
We now give some examples:

\begin{example}We work over the field of complex numbers. If $X=\ppp^1$ then
$\Pic(X) = \zzz$ and $H^0(G,\Pic(X)) =\zzz$ as well.
Since $\Pic(X)_0=0$, the answer to the question is
``yes'' if and only if there
is a fixed divisor of any given degree
(since there is only one class with a given degree).
We show that this can fail.
Fix an embedding of $G=A_5$ into $PGL(2,\ccc)$ 
(to get an action of $G$ on $\ppp^1$).
In this case, inertia groups (i.e. stabilizers of 
points on the curve) are always cyclic, so the only possibilities 
are of order $1$, $2$, $3$ and $5$. 
Thus the possible orbit sizes
are $12$, $20$, $30$ and $60$ and any $G$-fixed 
divisor has degree a multiple
of the greatest common divisor of these numbers (i.e. $2$).
Now Theorem~\ref{thrm:main} implies that
$H^2(A_5,\ccc^\times)=\Pic(X)^G/\Div(X)^G= \zzz/2\zzz$,
which of course is well-known (see e.g. \cite{K}, pp. 245--246).

Now let $F$ be an algebraically closed field of characteristic~$2$.
Then $H^2(A_5,F^{\times})=0$, 
since the $2'$ part of the Schur multiplier of $A_5$ is zero.
Fix an embedding of $A_5$ into $PGL(2,F)$ so that $A_5$ acts 
on $X=\ppp^1$.
It follows from Theorem~\ref{thrm:main} that 
$\Div(X)^G=\Pic(X)^G=\zzz.$ Indeed, the inertia groups are
$A_4$ and cyclic of order $5$, and the gcd of their indices is one.
\end{example}

\begin{example}Let $F$ be any closed field of characteristic not $2$.
Let $G$ be an elementary abelian group of order $4$ acting
on the curve $\ppp^1(F)$. One can argue as in the previous example
that $H^2(G,F^{\times})=\zzz/2\zzz$.
\end{example}

\begin{example}Let $F$ be any closed field of characteristic not $2$.
Let $G$ be an elementary abelian group of order $2^r$ acting
on the curve $X$. If $G$ acts fixed point freely on $X$, then
$B = 2^r\zzz$ in the notation above.  Since $A/B$ is cyclic
and embeds in $H^2(G,F^{\times})$, a group of exponent $2$,
it follows that every $G$-invariant divisor class on $X$
has degree a multiple  of $2^{r-1}$. We show below that
there may be such a class or there may not be.
Suppose that $G$ does not act fixed point freely on $X$.
Since $G$ acts tamely on $X$, all inertia groups are cyclic
and so have index a multiple of $2^{r-1}$, it follows
that $B=2^{r-1}\zzz$ and so every $G$-invariant divisor
class has degree a multiple of $2^{r-2}$.
\end{example}

We have already observed that
there are exact sequences:
\begin{equation}\label{eq1}
\Div^G(X) \to \Pic^G(X) \to H^2(G,F^\times) \to 0\rm {\ and}
\end{equation} 
\begin{equation}\label{eq2}
\Div_0^G(X) \to \Pic_0^G(X) \to H^2(G,F^\times) 
\to A/B \to 0,
\end{equation}
where $A =\deg(\Pic(X)^G)$ and $B=\deg(\Div(X)^G)$. Note that 
$B$ is generated by the $\gcd$ of $|G:\Stab_G(x)|$ as $x$ ranges over
$X$.

It is convenient to define $C=C(X,G)=A/B$. We summarize some
properties of this group.

\begin{lemma}\label{lem:C} We have:
\begin{enumerate}
\item The group $C=A/B$ is finite and cyclic.
\item $|C|$ divides the gcd of the $|G|/|I|$ as $I$ ranges over the
inertia groups.
\item  $C=1$ if there is a totally ramified point in $X$.
\end{enumerate}
\end{lemma}
Only the third statement deserves comment. Indeed, if the point $P\in X$ 
is totally ramified, then $[P]$ is a $G$-invariant degree one divisor.

David Saltman has shown us how to explicitly compute the last map 
$H^2(G,F^\times) \to A/B$ in the second statement.
He has kindly allowed us to include his observations here.
Let $n$ be an integer coprime to the characteristic of $F$.
Let $c$ be a $2$-cocycle in $H^2(G,\mu_n)$.
Then $c=\delta(d_g)$, where $d_g$ is a coboundary with values in 
$F(X)^\times$ by Tsen's theorem.
Since $d_g^n$ is trivial in $H^1$, it follows that 
\begin{equation}\label{one}
d_g^n=g\theta/\theta
\end{equation}
for some function $\theta$ on $X$.
Now $F(X)(\theta^{1/n})/F(X)^G$ is a Galois cover that
realizes $c$.
The class of $d_g$ vanishes in $H^1(G,\Div(X))$  and so we can write
(switching now from multiplicative to additive notation)
$$(d_g)= gD-D$$ 
for a divisor $D$. Now multiply by $n$ and use \eqref{one}
to get 
$ (g\theta) - (\theta) = ngD-D$. 
Hence 
$$(\theta) - D\mbox{\quad is a $G$-invariant divisor.}$$
\medskip
Now $D$ is not well-defined -- but its degree is well-defined modulo
$b$.
The map in question $H^2(G,F^\times) \to A/B$ sends 
(the class of) $c$ to the degree of $D$ via the construction above.
\medskip

Note that the question of how many  $G$-invariant divisor classes
are not in the image of a $G$-invariant divisor depends only on the group.
We show that this is not the case for the analogous question
for degree $0$ divisors and divisor classes. The answer depends on the
curve even if we assume that $X \rightarrow X/G$ is unramified.

We have already set $C(X,G)$ to be
$\deg(\Div(X)^G)/\deg(\Pic(X)^G)$.
By definition, $C(X,G)$ equals $A/B$.
We already observed that the term $H^2(G,F^\times)$ 
in $\eqref{eq1}$ does not depend on the curve $X$.
So it seems natural to ask whether the term $C=C(X,G)$ 
only depends on $G$ and the ramification data of $X$.
We next show Theorem~\ref{klein} that indeed, 
this fails over the complex field even
for $G$ equals the Klein four group acting without fixed points.

We will need two lemmas to prove the theorem.

\begin{lemma}\label{godown}
Let the finite group $G$ act on the curve $X$.
Let $N$ be a normal subgroup of $G$ that acts on $X$ without fixed
points. Let $Q$ be the quotient group
and let $Y$ be the quotient curve. 
Then $|C(Y,Q)|$ divides $|C(X,G)|$.
\end{lemma}
\begin{proof} A $Q$-invariant divisor class on $Y$ of degree $d$ lifts to
a $G$-invariant divisor class of degree $|N|d$ on $X$. However, 
for divisors, you can also go down: a 
$G$-invariant divisor of degree $|N|d$ on $X$ always
comes from a $Q$-invariant divisor on $Y$ of degree $d$.
\end{proof}

\begin{lemma}  Let $X$ be a curve over the algebraically closed
field $F$ of characteristic not $2$.  Let $G$ be an elementary abelian
$2$-group of order $2^r$ with $r \ge 3$.  Then
the cokernel of $\Div(X)^G \rightarrow \Pic(X)^G$ is an 
elementary abelian group of order $2^m$ where $m=r(r-1)/2$.
\end{lemma}
\begin{proof}  The lemma follows from the fact that $H^2(G,F^\times)$ equals
$\wedge^2(G)$ and from Theorem~\ref{thrm:main}.
\end{proof}
\begin{theorem}\label{klein} Let $K$ be the Klein four group. 
\begin{enumerate}
\item There exists a curve
$E$ over the field of complex numbers such that $K$ acts
without fixed points on $E$ and $C(E,K)=\zzz/2\zzz$.
\item There exists a curve
$Y$ over the field of complex numbers such that $K$ acts
without fixed points on $Y$ and $C(Y,K)=1$.
\end{enumerate}
\end{theorem}
\begin{proof}
Let $E$ be an elliptic curve over $\ccc$ 
and let $K$ act by translations on $Y$.
Note that indeed $K$ acts without fixed points.
In particular, any $K$-invariant divisor on $E$ has degree $4$.

Let $\{0,P,Q,R,\}$ be the of two-torsion of $E(\ccc)$.
Using the group law on $E$ we have $0+P= Q+R$, which implies
that as divisors  $(0)+(P)\sim (Q)+(R).$ Note further
$E(\ccc)$ and in particular $K$ acts trivially on $\Pic_0(E).$ 
Thus there is an invariant divisor class of degree $2$.
Thus $C(E,K)=\zzz/2\zzz$ and the first part is proved.

Let $r \ge 3$ be a positive integer.  Let $J$ be the affine
group $AGL(r,2)$. The subgroup of translations $G$ is elementary abelian
of order $2^r$. Note that $G$ is normal in $J$. 
Now $J$ acts on some curve $X$ with
the cover $X \rightarrow X/J$ unramified (in characteristic
$0$, this follows from the description of the fundamental
group -- in positive characteristic, this follows from
\cite{stevenson}).
Now we note that it follows from our proofs above
that the map from
$H^1(G,\Prin(X)) \rightarrow H^1(G,\Div_0(X))$ is $J$-equivariant.
Since $H^1(G,\Prin(X)) \cong \wedge^2(G)$ is $J$-irreducible
and not cyclic (since $r \ge 3$), it 
follows that $C(X,G)$ is trivial.
Now let $N$ be a subgroup of $G$ of index $4.$ We have $K=G/N$ and
set $Z=X/N$. Now lemma \ref{godown} implies that  $C(Z, K)$ is
trivial. So we have constructed two curves with the same group acting
without fixed points such that the cokernels 
of $\Div_0(X)^G \rightarrow \Pic_0(X)^G$ are different.
The theorem is proved.
\end{proof}

One can choose the two curves 
to have the same genus as well.

\section{Varieties}
In this section, $X$ is a non-singular, projective variety 
over the field $F$ of complex numbers.
Let $D$ and $D'$ be $(n-r)$-cycles on $X$.
Recall that if $D$ and $D'$ are rationally equivalent then they are
algebraically equivalent, if $D$ and $D'$ are algebraically 
equivalent then they are homologically equivalent,
if $D$ and $D'$ are homologically equivalent
then they are numerically equivalent
(\cite{Har}, \S V.1, Exer. 1.7, and \cite{Ful1}, \S 19.3).
Using the notation from the introduction above, 
for $D\in \Div(X)$,
\[
\begin{split}
&[D]_r=D+Rat^r(X)\subset [D]_a=D+\Alg^r(X)\subset \\
&[D]_h=D+\Hom^r(X)
\subset [D]_n=D+\Num^r(X).
\end{split}
\]
When specifying which notion of equivalence is not needed,
we drop the subscript. 
\begin{lemma}
\label{lemma:equivar}
If $[D]_r$ is $G$-equivariant then so are
$[D]_a$, $[D]_h$, and $[D]_n$.
\end{lemma}
\begin{proof}
If $[D]$ is $G$-equivariant if and only if for all $g\in G$,
$gD-D$ is equivalent to $0$. 
Since $Rat^r(X)\subset \Alg^r(X)\subset \Hom^r(X)\subset \Num^r(X)$,
we have the result claimed in the lemma.
\end{proof}
Consider again the exact sequences (\ref{eqn:longer}) and
(\ref{eqn:long}), where $X$ is now a variety. 
\subsection{Divisors}
We start with the analog of a result in the
previous section.
\begin{lemma}
The map
$H^1(G,\Prin(X)) \rightarrow H^2(G,F^\times)$
is injective.
\end{lemma}
\begin{proof}
By Hilbert's Theorem 90, $H^1(G,F(X)^\times )=1$, 
so this follows from (\ref{eqn:longer}).
\end{proof}
As in the case of curves, it follows that there is
a $G$-equivariant representative
of a $G$-equivariant divisor class if 
the Schur multiplier of $G$ is trivial.
\begin{theorem}  
\label{thrm:main2}
If $H^2(G,F^\times)=1$ then the map
$\Div(X)^G \rightarrow \Pic(X)^G$
is surjective. 
In other words, if $H^2(G,F^\times)=1$ and $[D]_r$ is $G$-invariant
then there always a $D'\in [D]_r$ which is also $G$-equivariant.
\end{theorem}

\begin{corollary} 
\label{cor:11}
If every $\ell$-Sylow subgroup of $G$ is cyclic 
(for every prime $\ell$ dividing $|G|$)
then for each $G$-invariant divisor class $[D]_r$
there is always a $D'\in [D]_r$ which is $G$-equivariant.
\end{corollary}
Recall
there is a natural 1-1 correspondence between invertible 
subsheafs ${\mathcal L}(D)$ of the sheaf ${\mathcal K}$ of total quotient rings on
$X$ and divisors $D$ (\cite{Har}, \S II.6). 
\begin{corollary}
If the hypothesis to Corollary \ref{cor:11} holds then for each
$G$-equivariant invertible subsheaf ${\mathcal L}$ of ${\mathcal K}$ 
on $X$ there is a $G$-equivariant divisor
$D$ of $X$ such that ${\mathcal L}={\mathcal L}(D)$.
\end{corollary}
\begin{proof}
Use the above correspondence and the previous
corollary.
\end{proof}
The long exact sequence (\ref{eqn:longer})
does not help to determine the image
of the map
\[
\Div(X)^G\rightarrow (\Div(X)/\Alg^1(X))^G,
\]
or of the map
\[
\Div(X)^G\rightarrow (\Div(X)/\Num^1(X))^G.
\]
However, in some special cases, one can say more.
\begin{theorem} Let $X$ be a $K3$ surface.
If $H^2(G,F^\times)=1$ then the map
$\Div(X)^G \rightarrow (\Div(X)/\Alg^1(X))^G$,
is surjective. 
In other words, 
if $H^2(G,F^\times)=1$ and $[D]_a$ is $G$-invariant
then there always a $D'\in [D]_a$ which is $G$-equivariant.
The analogous statements with $\Alg^1$ replaced by
$\Hom^1$ or $\Num^1$ also hold.
\end{theorem}
\begin{proof} 
For $K3$ surfaces, $[D]_r= [D]_a= [D]_h=[D]_n$
(\cite{SD}, \S 2.3),
so this is a consequence of the analogous result 
proven previously for rational equivalence.
\end{proof}
\subsection{Toric varieties}
In this subsection, we show that the answer to the 
question asked in our introduction 
is ``yes'' for rational cycle classes on
toric varieties, at least if $G$ consists only of
``toric automorphisms''. In the toric case, the 
Chow groups can be described in terms of the
combinatorial geometry of the fan associated to the
variety. 
Let $\Delta$ be a finite, non-singular, strongly polytopal, complete 
fan associated to an
integral lattice $L$ in $\rrr^n$ and let $X=X(\Delta)$ denote the 
associated
toric variety. (See Definition V.4.3 of \cite{ewa} for the 
definition of strongly polytopal.)
The variety $X$ is complete and projective and
contains a torus $T$ densely.
A {\bf toric automorphism} is a $T$-equivariant automorphism of
$X(\Delta)$. 
\begin{theorem}
Suppose that 
\begin{enumerate}
\item[(a)] $G\subset {\rm Aut}(X)$
is a finite subgroup of the group of toric
automorphisms of $X$;  and 
\item[ (b)] the 
class $[Z]\in Z^r(X)/Rat^r(X)$ is $G$-equivariant.
Then there is a $Z'\in [Z]$ which is $G$-equivariant.
\end{enumerate}
\end{theorem}
\begin{proof}
Toric automorphisms of $X$ correspond to automorphisms
of the lattices which also preserve $\Delta$ (Theorem 1.13 in Oda \cite{O}).
We know by the Proposition in \S 5.1 in \cite{Ful2} that
$[Z]$ is an integral combination of the classes of the orbit closures
$V(\sigma)$, $\sigma$ a cone in $\Delta$:
\[
[Z]_r=\sum_i n_i [V(\sigma_i)]_r,\ \ \ \ \ \ n_i\in \zzz.
\]
But such a combination of classes is fixed by $G$ only if
only if the cones $\{\sigma_i\}$ decomposes into a disjoint
union of orbits and the $n_i$'s are constant on each
of these primitive $G$-orbits.
In this case, the cycle $\sum_i n_i V(\sigma_i)$ itself is 
$G$-equivariant.
\end{proof}

\begin{remark}
 Let $X$ be a non-singular variety for which
$H^2(X,\zzz)$ is torsion-free. (This
is true for toric surfaces, for example \cite{Ful2}, \S 3.4.) 
Then algebraic, homological and numerical
equivalence of divisors agree (\cite{Ful1}, \S 19.3.1).
Consequently, if
$Z^r(X)^G \rightarrow (Z^r(X)/\Alg^r(X))^G$
is surjective then so is
$Z^r(X)^G \rightarrow (Z^r(X)/\Num^r(X))^G$,
and conversely.
\end{remark}


\begin{thebibliography}{99}
\bibitem[Ewa]{ewa}
G. Ewald, {\bf Combinatorial complexity and algebraic geometry},
Springer-Verlag, 1996.
\bibitem[Ful1]{Ful1} W. Fulton, {\bf Intersection theory},
2$^{nd}$ ed., Springer-Verlag, 1998.
\bibitem[Ful2]{Ful2} W. Fulton, {\bf Introduction to
toric varieties}, Princeton Univ. Press, 1993.
\bibitem[H]{H} M. Hall, {\bf Theory of groups}, AMS Chelsea,  
1976.
\bibitem[Har]{Har} R. Hartshorne, {\bf Algebraic geometry},
Springer-Verlag, 1977.
\bibitem[K]{K} G. Karpilovsky, {\bf The Schur multiplier},
Oxford Univ. Press, 1987.
\bibitem[Lo]{Lo} K. Lonsted, ``On $G$-linebundles and $K_G(X)$,''
J. Math. Kyoto Univ. \underbar{23}(1983)775-793.
\bibitem[O]{O} T. Oda, {\bf Convex bodies and algebraic geometry},
Springer-Verlag, 1985.
\bibitem[SD]{SD} 
B. Saint-Donat, ``Projective models of K-3 surfaces,'' 
Amer. J. of Math., vol. \underbar{96} (1974) 602-639.
\bibitem[S]{S} J.-P. Serre, {\bf Local fields}, Springer-Verlag, 
1979.
\bibitem[Sh]{Sh} S. Shatz, {\bf Profinite groups, arithmetic,
and geometry}, Princeton Univ. Press, 1972.
\bibitem[Ste]{stevenson}  K. Stevenson,  
Galois groups of unramified covers of projective curves in characteristic $p$,
J. Algebra 182 (1996),  770--804.
\bibitem[Sti]{St} H. Stichtenoth, {\bf Algebraic function
fields and codes}, Springer-Verlag, 1993.
\end{thebibliography}
\end{document}